\documentclass[11pt]{amsart}

\newlength{\basicwidth} 

\newlength{\shortbasicwidth} 

\newlength{\basicheight} 

\begin{document}

\numberwithin{equation}{section}

\begin{center}
\bf{\Large{A Hardy-Littlewood Integral Inequality on  Finite Intervals}} \\
\end{center}

\begin{center}
\bf{\Large{with a Concave Weight}}
\end{center}

\par\vspace*{2\baselineskip}\par
\centerline{\large  ({\em To appear in: Periodica Mathematica Hungarica})}

\vspace{1cm}
\begin{center}
HORST ALZER$^a$ \, and \, MAN KAM KWONG$^b$\footnote{The research of this author is fully supported by 
a grant from the Research Grants Council
of the Hong Kong Special Administrative Region, China (Project No.\ PolyU
5012{/}10P).}
\end{center}

\vspace{0.2cm}
\begin{center}
${}^a$ \emph{Morsbacher Str. 10, 51545 Waldbr\"ol, Germany}\\
\tt{H.Alzer@gmx.de}
\end{center}

\vspace{-0.2cm}
\begin{center}
${}^b$ \emph{Department of Applied Mathematics, The Hong Kong Polytechnic University, \\
Hunghom, Hong Kong}\\
\tt{mankwong@polyu.edu.hk}
\end{center}

\vspace{2cm}
{\bf{Abstract.}}
We prove: For all concave functions $ w: [a,b] \rightarrow [0,\infty)$ and for all functions $f \in C^2[a,b]$ with $f(a)=f(b)=0$ we have
\[
  \left(  \int_{a}^{b} w(x) f'(x)^2 \,dx \right) ^2 
\leq  
\Bigl(  \int_{a}^{b} w(x) f(x)^2 \,dx \Bigr)   \, \Bigl( \int_{a}^{b} w(x)f''(x)^2 \,dx\Bigr).
\]
Moreover, we determine all  cases of equality.

\vspace{1cm}


\newpage

\setcounter{equation}{0}\section{Introduction}

In  1932, Hardy and Littlewood  \cite{HL} established integral inequalities
of the form
\begin{equation}  \left(  \int_J f'(x)^2 \, dx \right) ^2 \leq  k(J) \left(  \int_J f(x)^2 \,d x \right) \left( \int_J f''(x)^2  \,d x \right) ,    \end{equation}
where $ J $ is either $ \mathbf R=(-\infty ,\infty ) $ or $ \mathbf R^+=(0,\infty ) $, and
$ f $ is a twice-differentiable function such that $ f,f''\in L^2(J) $, with
the best possible constants
\begin{equation}  k(\mathbf R) = 1, \quad  k(\mathbf R^+) = 2 .   \end{equation}
The square integrability of $ f' $ and
the existence of a finite constant $ k(J) $ independent of $ f $
are implicitly parts of the conclusion. 

The significance of the result is in providing an estimate of the ``size''
of the derivative of a function when bounds on the ``sizes'' of the function and
its second derivative are known, with potential applications in the
study of differential equations.

The above results of Hardy and
Littlewood are reproduced in Sections 7.9 and 7.8 of the classical
work ``Inequalities'' by Hardy, Littlewood, and P\'olya \cite{HLP}.
The case $ J=\mathbf R $ is an immediate consequence of integration
by parts and the Cauchy-Schwarz inequality. The 
$ J=\mathbf R^+ $ case has a lengthy proof using calculus of variations.
Shorter proofs are found later; see, for example, \cite{KZ5}.

Inequalities of the form (1.1)
 are special cases of a more general result, called the Landau
inequalities, obtained by replacing the $ L^2 $ norms (of $ f' $ on the lefthand
side and of $ f $ and $ f'' $ on the righthand side) by other 
$ L^p $ norms, with $ 1\leq p<\infty  $. We note that now the best
constant $ k(p,J) $ depends both on $ J $ and $ p $. It is known
that $ k(p,\mathbf R)<k(p,\mathbf R^+) $ for all $ p $.
The exact values of $ k(p,J) $
are known only for $ p=1 $ (due to Berdyshev \cite{B})
and $ p=\infty  $ (due to  Hadamard and  Landau, respectively), 
besides the Hardy-Littlewood result (1.2) for $ p=2 $.
We only mention here that
$ k(1,\mathbf R) = 2  $.
For details on the other results and further readings, 
we refer the readers to the monograph \cite{KZ4}.

Many further extensions have been obtained.
 Gabushin used three (possibly) different norms,
$ L^p,L^q $ and $ L^r $ for $ f,f' $ and $ f'' $, respectively, and 
also studied higher derivative analogues.
Kato \cite{K} replaced $ L^2(J) $ with a Hilbert space, 
and the derivatives of $ f $ with $ Af $ and $ A^2f $, respectively,
where $ A $ is an m-dissipative operator.
Everitt \cite{Ev} substituted the differentiation operator with a general
self-adjoint second-order differential operator, and showed that the 
corresponding best constant
can be determined using the $m_\lambda $ function that 
Titchmarsh introduced in the Weyl
theory of limit-point{/}limit-circle classification of second-order
differential operators. The resulting inequality has been dubbed the
HELP (Hardy-Everitt-Littlewood-P\'olya) inequality. 
See \cite{EE} for a survey of this topic.

In this paper we are concerned with a weighted form of (1.1)
which holds for functions defined on finite intervals and subjected to suitable
boundary conditions. More precisely, we determine the best possible constant
$\kappa$  such that we have for all concave functions 
$ w: [a,b]\rightarrow [0,\infty)$ and for all $f\in C^2[a,b]$ with
$f(a)=f(b)=0$: 
\begin{equation}
  \left(  \int_{a}^{b} w(x) f'(x)^2 \,dx \right) ^2 
\leq \kappa \cdot \,   
\Bigl(  \int_{a}^{b} w(x) f(x)^2 \,dx \Bigr)   \, \Bigl( \int_{a}^{b} w(x)f''(x)^2 \,dx\Bigr).
\end{equation}

Two previous works of the second author pertain to the current paper.
First, we note that if $ J $ is replaced by a finite interval $ [a,b] $, an
inequality of the form (1.1), as well as its more general Landau 
analogue,
cannot hold, as shown by the counterexample
$ f(x)=x $. However, if one of the three sets of
boundary conditions
$$  
f(a) = f(b) = 0,   
\quad{ f(a) = f'(b) = 0,}    
\quad\mbox{or}\quad{  f'(a) = f(b) = 0,}
$$    
is imposed, then (1.1)
(or its Landau analogue) is preserved with 
the same best constant for $ J=\mathbf R $.
Indeed, it was shown by Kwong and Zettl \cite{KZ2}
that any one of these finite-interval results
is equivalent to the original result for $ J=\mathbf R $.
The case $ J=\mathbf R^+ $ also has a finite-interval equivalent associated with the
boundary condition $f'(a)/f(a) = f'(b)/f(b)$. For details, see \cite{KZ4}.

In \cite{KZ5}, it was proved that a Landau-type inequality holds 
when a weight function $ w $ is added in defining the norms of functions, 
provided that $ w $ is an increasing function:
\begin{equation}  \|f'\|^2 \leq  k(p,\mathbf R^+) \, \|f\|\cdot \|f''\|,    \end{equation}
where
$$
 \|f\| = \left(  \int_J  w(x)|f(x)|^p \, dx \right) ^{1/p}.   
$$
We are not claiming that the constant $ k(p,\mathbf R^+) $ that appears on the 
righthand side of (1.4)
is the best possible constant for a given weight $ w$; it is, however,
best possible over all choices of increasing $ w $.
Another important point to note is that, unlike the classical Landau 
inequalities, 
the same constant is required in (1.4),
irrespective of whether
$ J=\mathbf R $ or $ \mathbf R^+ $.

\section{Lemmas}

In this section, we collect a few basic facts on concave functions, which we 
need for our purposes. Proofs for the first three lemmas given below as well as
  more information on this subject can be found, for example, in the monographs
Niculescu and Persson \cite{NP}, Roberts and Varberg \cite{RV}, and Royden and Fitzpatrick
\cite[Section 6.6]{RF}.

We recall that a function $f: I \rightarrow \mathbf{R}$, where $I\subset \mathbf{R}$ is an
interval, is said to be concave if
$$
\lambda f(x) +(1-\lambda) f(y)\leq f(\lambda x+(1-\lambda)y)
$$
for all $x,y \in I$ and $\lambda \in (0,1)$. If $-f$ is concave, then $f$ 
is called convex.

{\bf{Lemma 1.}} \emph{If $f: [a,b]\rightarrow \mathbf{R}$ is concave, then $f$ is continuous on $(a,b)$ and the limits }
$$
\lim_{x\to a+} f(x) \quad{and} \quad{\lim_{x\to b-} f(x)}
$$
\emph{exist.}

From now on, we always assume that (if necessary) $f$ has been modified
such that $f(a)=\lim_{x\to a+} f(x)$ and  $f(b)=\lim_{x\to b-} f(x)$ exist,
so that $f$ is continuous on $[a,b]$.

The left and right derivatives of a function are defined by
$$
f_{-}'(c)= \lim_{x\to c-}\frac{f(x)-f(c)}{x-c},
\quad{f_{+}'(c)= \lim_{x\to c+}\frac{f(x)-f(c)}{x-c}}.
$$

{\bf{Lemma 2.}} \emph{If $f: [a,b]\rightarrow \mathbf{R}$ is concave, then 
$f_{-}'$ and $f_{+}'$ exist and are decreasing on $(a,b)$.}

{\bf{Lemma 3.}} \emph{If $f: (a,b)\rightarrow \mathbf{R}$ is concave, then 
we have for $c,x\in (a,b)$:}
$$
f(x)-f(c)=\int_c^x f_{-}'(t)dt=\int_c^x f_{+}'(t)dt.
$$

The next result is nothing more than integration by parts for
integrable functions. It is given, for example, in \cite[p. 32]{T} and in
\cite[p. 128, Problem 52]{RF}.

{\bf{Lemma 4.}} \emph{Let  $w: [a,b] \rightarrow \mathbf{R}$ be concave.
If $f \in C^1[a,b]$ with $f(a)=f(b)=0$, then}
$$
\int_a^b w(t) f'(t)dt=-\int_a^b w_{-}'(t) f(t)dt.
$$

{\bf{Lemma 5.}} \emph{Let $w: [a,b]\rightarrow [0,\infty)$ be  concave. 
Then there exists a sequence of non-negative concave functions 
$(w_n)$ $(n=1,2,...)$ with $w_n\in C^2[a,b]$ such that $w_n$ converges 
uniformly to $w$  on $[a,b]$.   }

\begin{proof}
A result of Bremermann \cite{Br} states that if $f$ is convex on $[a,b]$, then
for every $k\in\mathbf{N}$ there exists a sequence of $k$-times differentiable convex functions $(f_n)$ with $f_n(x)\geq f_{n+1}(x)\geq f(x)$ and 
$f_n$ converges to $f$ on $[a,b]$. Applying Dini's theorem we conclude
that the convergence is uniform.  This implies that there
exists a sequence of concave functions $(w_n)$ with $w_n\in C^2[a,b]$
such that $w_n$ converges uniformly to $w$ on $[a,b]$. 
If it happens that $w_n\geq 0$ on $[a,b]$ for all $ n $, then 
$  (w_n)  $ is our
desired sequence.

In the contrary case, we modify $ w_n $ as follows.
First we note that,
for every  $\epsilon >0$  there is a  natural number $n_0$ such that for all $n \geq n_0$ and for all $x\in[a,b]$:   
\begin{equation}
|w(x)-w_n(x)| < \epsilon/2.   \label{wn}
\end{equation}
We assume that $w_n$ is negative somewhere in $[a,b]$ and set
$$
w_n(x_0)=\min_{a\leq x\leq b} w_n(x).
$$
Then, $w_n(x_0)<0$. From (2.1) we obtain
$$
w(x_0)-w_n(x_0)<{\epsilon}/{2}.
$$
Since $w(x_0)\geq 0$, it follows that
$$
w_n(x_0)>w(x_0)-{\epsilon}/{2}\geq - {\epsilon}/{2}.
$$
Let $ w^*_n(x)=w_n(x)-w_n(x_0)$. 
Then, $w_n^*$ is non-negative and concave on $[a,b]$.
Moreover, we get for $n\geq n_0$ and $x\in[a,b]$:
$$
|w(x)-w_n^*(x)|\leq |w(x)-w_n(x)|+|w_n(x_0)|<
{\epsilon}/{2}+{\epsilon}/{2}={\epsilon}.
$$
This implies that $w_n^*$ converges uniformly to $w$ on $[a,b]$.
\end{proof}

\section{Main Result}

We are now in a position to present our main result. The following theorem 
reveals that the best possible constant factor in (1.3) is given by $\kappa=1$.

{\bf{Theorem.}}
\emph{For all concave functions $ w: [a,b] \rightarrow [0,\infty)$ and for all functions $f \in C^2[a,b]$ with $f(a)=f(b)=0$ we have}
\begin{equation}
  \left(  \int_{a}^{b} w(x) f'(x)^2 \,dx \right) ^2 
\leq  
\Bigl(  \int_{a}^{b} w(x) f(x)^2 \,dx \Bigr)   \, \Bigl( \int_{a}^{b} w(x)f''(x)^2 \,dx\Bigr).
\end{equation}
\emph{Let $w\not \equiv 0$ and  $f\not \equiv 0$. Then, the sign of equality 
holds in} (3.1) \emph{if and only if} 
\begin{equation}  
f(x)=\lambda \sin \left(  \frac{n\pi (x-a)}{b-a} \right) , 
\end{equation}
\emph{where $\lambda$ is a real number and $n$ is a natural number,
and  $w$ is linear in each of the subintervals}
\begin{equation}
J_k=\Bigl[a+\frac{(k-1)(b-a)}{n}, a+\frac{k(b-a)}{n}   \Bigr],
\quad{k=1,...,n}. 
\end{equation}

\begin{proof}
From  Lemma 5, we conclude that it  
suffices to  prove  inequality (3.1) under the additional assumption that the 
weight function $ w \in  C^2[a,b]$.

Integration by parts gives
\begin{equation}
\int_a^b w''(x) f(x)^2 dx = \Bigl[w'(x)f(x)^2\Bigr]_a^b 
-2 \int_a^b w'(x) f(x) f'(x)^2 dx   
\end{equation}
and
\begin{equation}
\int_a^b w'(x) f(x) f'(x) dx = \Bigl[w(x)f(x) f'(x)\Bigr]_a^b 
- \int_a^b w(x)[ f(x) f''(x)+f'(x)^2] dx.   
\end{equation}
Since $f(a)=f(b)=0$, we obtain from (3.4) and (3.5):
\begin{equation}
\int_a^b w''(x) f(x)^2  dx =  
2 \int_a^b w(x)[ f(x) f''(x)+f'(x)^2] dx.   
\end{equation}
Using  $-w''(x)\geq 0$ and  (3.6), we find
\begin{equation}
\int_a^b w(x)  f'(x)^2 dx \leq 
 \int_a^b w(x) f'(x)^2 dx-\frac{1}{2}
\int_a^b w''(x) f(x)^2  dx   
=
\int_a^b \bigl( \sqrt{w(x)} f(x) \bigr) \cdot 
\bigl(- \sqrt{  w(x)}  f''(x) \bigr) dx.  
\end{equation}
The Cauchy-Schwarz inequality yields
\begin{equation}
\int_a^b \bigl( \sqrt{w(x)} f(x) \bigr) \cdot 
\bigl(- \sqrt{  w(x)}  f''(x) \bigr) dx
\leq
\Bigl(\int_a^b w(x) f(x)^2 dx  \Bigr)^{1/2}
\Bigl(\int_a^b w(x) f''(x)^2 dx  \Bigr)^{1/2}.   
\end{equation}
Combining (3.7)   and (3.8)  leads to (3.1).

\vspace{0.1cm} 
Next, we discuss the case of equality.
A short calculation reveals that equality holds in (3.1) if
$ f $ is given by (3.2) and 
if $w$ is linear in each of the subintervals $J_k$ given in (3.3).

Now, we assume that the sign of equality is valid in (3.1). Let 
$ w \not \equiv 0$ and
$ f \not \equiv 0$. Applying  Lemma 5, (3.7), and (3.8) yields
$$
\int_a^b w_n(x)  f'(x)^2 dx
\leq
- \int_a^b w_n(x) f(x) f''(x) dx
\leq
\Bigl(\int_a^b w_n(x) f(x)^2 dx  \Bigr)^{1/2}
\Bigl(\int_a^b w_n(x) f''(x)^2 dx  \Bigr)^{1/2}.
$$
Letting $n\rightarrow \infty$ gives
\begin{equation}
\int_a^b w(x)  f'(x)^2 dx
\leq
- \int_a^b w(x) f(x) f''(x) dx
\leq
\Bigl(\int_a^b w(x) f(x)^2 dx  \Bigr)^{1/2}
\Bigl(\int_a^b w(x) f''(x)^2 dx  \Bigr)^{1/2}.
\end{equation}
By assumption, the expressions on the left-hand side and right-hand side
are equal. This implies that the sign of 
equality holds in the second inequality of (3.9). It follows that
$$
\sqrt{w(x)} f(x)=\lambda_0 \sqrt{w(x)} f''(x)
$$
for some constant $\lambda_0$. Since $w$ is non-negative and concave with 
$w \not\equiv 0$, we conclude that $w$ is positive on $(a,b)$. Hence,
$$
f(x)=\lambda_0 f''(x)
$$
with $\lambda_0 \neq 0$. Together with the boundary conditions $f(a)=f(b)=0$,
we see that $f$ must be of the form given in (3.2).

Since $f$ is equal to $0$
at the end-points of each subinterval $J_k$, we conclude
that the hypothesis of the Theorem  is satisfied for each $k$. Hence,
 for all $k\in \{1,2,...,n\}$ we have   
\begin{equation}
  \left(  \int_{J_k} w(x) f'(x)^2 \,dx \right) ^2 
\leq  
\Bigl(  \int_{J_k} w(x) f(x)^2 \,dx \Bigr)   \, \Bigl( \int_{J_k} w(x)f''(x)^2 \,dx\Bigr).
\end{equation}
We claim that equality must hold in (3.10) for each $k$.

The following equivalence is easy to verify:

\begin{quote}\em
If $A$, $B$, and $C$ are positive numbers, then 
\begin{equation}
B^2 \leq AC \quad{\Longleftrightarrow} \quad{B \leq \epsilon A +\frac{1}{4\epsilon }C}
\quad\mbox{for all} \quad{\epsilon >0}.
\end{equation}
Moreover, ``$<$'' holds on the left-hand side if and only if it holds
on the right-hand side. 
\end{quote}

Applying (3.11) 
reveals that (3.10) is equivalent to
\begin{equation}
    \int_{J_k} w(x) f'(x)^2 \,dx 
\leq  
\epsilon   \int_{J_k} w(x) f(x)^2 \,dx + \frac{1}{4\epsilon}  \int_{J_k} w(x)f''(x)^2 \,dx \quad\mbox{for all} \quad{\epsilon >0}.
\end{equation}
If strict inequality in (3.10) holds for at least one $k$, then the 
corresponding inequality (3.12) must also be strict.
After adding up (3.12) for all $k$, we obtain
$$
    \int_{a}^{b} w(x) f'(x)^2 \,dx 
<  
\epsilon   \int_{a}^{b} w(x) f(x)^2 \,dx + \frac{1}{4\epsilon}  
\int_{a}^{b} w(x)f''(x)^2 \,dx \quad\mbox{for all} \quad{\epsilon >0}.
$$
This implies that the sign of equality does not hold in (3.1), a contradiction.
 
It remains to show that $w$ is linear in each $J_k$ under the assumptions
that $f$ has the representation  (3.2) and that 
 equality is valid in (3.10) for each $k$. 
We insert (3.2) in (3.10) (with ``$=$'' instead of ``$\leq $'')
 and substitute
$$
x=a+\frac{(k-1)(b-a)}{n}+\frac{b-a}{n\pi}t.
$$
Then we get
\begin{equation}
\int_0^\pi W(t) \cos^2(t)dt=\int_0^\pi W(t) \sin^2(t)dt
\end{equation}
with
$$
W(t)=w\Bigl( a+\frac{(k-1)(b-a)}{n}+\frac{b-a}{n\pi}t  \Bigr).
$$
The function $W$ is concave.
We show that $W$ is linear on $[0,\pi]$, which implies that $w$ 
is linear on $J_k$.

The identity 
 $\cos^2(A)-\sin^2(A)=\cos(2A)$ and  (3.13) imply
\begin{equation}
\int_0^\pi W(t) \cos(2t)dt=0.
\end{equation}
Using Lemma 4 and (3.14) gives
\begin{equation}
\int_0^\pi W_{-}'(t) \sin(2t)dt=
\int_0^{\pi/2} W_{-}'(t) \sin(2t)dt+\int_{\pi/2}^\pi W_{-}'(t) \sin(2t)dt
=0.
\end{equation}  
We have
\begin{equation}
\int_{\pi/2}^\pi W_{-}'(t) \sin(2t)dt=-
\int_{0}^{\pi/2} W_{-}'(\pi-t) \sin(2t)dt,
\end{equation}
so that  (3.15) and (3.16) lead to
\begin{equation}
\int_0^{\pi/2} \Delta(t) \sin(2t)dt=0,
\end{equation}
where
$$
\Delta(t)=W_{-}'(t)-W_{-}'(\pi-t).
$$
Since $W_{-}'$ is decreasing on $(0,\pi)$, it follows that $\Delta \geq 0$ on $(0,\pi/2)$. Moreover, the function $t \mapsto -W_{-}'(\pi-t)$ is decreasing
on $(0,\pi/2)$, which implies that $\Delta$ is also decreasing
on $(0,\pi/2)$. Since $\sin(2t)>0$ for $t\in (0,\pi/2)$, we conclude
from (3.17) that $\Delta(t)=0$ for all $t\in (0,\pi/2)$. Hence,
\begin{equation}
W_{-}'(t)=W_{-}'(\pi-t) \quad{\mbox{for all }} \quad{t\in (0,\pi/2].}
\end{equation}
Let $\epsilon$ be any small positive number. Since 
$W_{-}'$ is decreasing we conclude from (3.18) that $W_{-}'$ 
must be a constant on $(\epsilon, \pi-\epsilon)$. This is valid
for any small $\epsilon >0$, which reveals that $W_{-}'$
is a  constant on $(0,\pi)$. Applying Lemma 3 reveals
 that $W$ is linear. 
\end{proof}

{\bf{Corollary.}}
\emph{For all concave and increasing functions $ w: [a,b] \rightarrow [0,\infty)$ and for all functions $f \in C^2[a,b]$ with $f(a)=f'(b)=0$ we have
\begin{equation}
  \left(  \int_{a}^{b} w(x) f'(x)^2 \,dx \right) ^2 
\leq  
\Bigl(  \int_{a}^{b} w(x) f(x)^2 \,dx \Bigr)   \, \Bigl( \int_{a}^{b} w(x)f''(x)^2 \,dx\Bigr).
\end{equation}
Let $w\not \equiv 0$ and  $f\not \equiv 0$. Then, the sign of equality 
holds in} (3.19) \emph{if and only if 
$$
f(x)=\lambda \sin \left(  \frac{(2n-1)\pi (x-a)}{2(b-a)} \right) , 
$$
where $\lambda$ is a real number and $n$ is a natural number,
and  $w$ is linear in each of the subintervals
$$
I_k=\Bigl[a+\frac{2(k-1)(b-a)}{2n-1}, a+\frac{2k(b-a)}{2n-1}   \Bigr],
\quad{k=1,...,n-1} 
$$
and is constant in the subinterval}
$$
I=\Bigl[a+\frac{2(n-1)(b-a)}{2n-1}, b   \Bigr].
$$

\begin{proof}
We extend $ w $ and $ f $ to $ [a,2b-a] $ so that each function is even with respect to $ x=b $, that is,
$$
w(b+\sigma ) = w(b-\sigma ) \quad{\mbox{and}}
\quad{   f(b+\sigma ) = f(b-\sigma )}  \quad  \mbox{for } \sigma \in(0,b-a].  
$$
Then,  $ w $ is non-negative and concave on $[a,2b-a]$
and $ f \in C^2 [a,2b-a] $ with $f(a)=f(2b-a)=0$.
So, $w$ and $f$
 satisfy the hypotheses of the Theorem on the extended domain $[a,2b-a]$
and the
conclusion follows from the Theorem.
\end{proof}

An obvious analogous result holds for functions $ f $ that satisfies $ f'(a)=f(b)=0 $.

\section{Concluding Remarks}

 In view of the Theorem  two questions arise naturally.

(i) Does the Theorem  remain valid if we remove the ``concavity'' requirement?

(ii)  Is there a corresponding result for convex weight functions? More 
precisely, given a convex function $w:[a,b]\rightarrow [0,\infty)$,
does there exist a constant $C$ such that  for all 
$f\in C^2[a,b]$ with $f(a)=f(b)=0$:
$$
\kappa(w,f)= \frac{ \left(  \int_{a}^{b} w(x) f'(x)^2 \,dx \right) ^2 }
{\Bigl(  \int_{a}^{b} w(x) f(x)^2 \,dx \Bigr)   \, \Bigl( \int_{a}^{b} w(x)f''(x)^2 \,dx\Bigr)}\leq C \, ?
$$

An example, constructed with the help of the MAPLE software,
 reveals that in both cases the answer is ``no''.

Let $ [a,b]=[0,1] $, $ w(x)=x^4 $, and  $\delta \in (0,1/2)$.
The  piecewise polynomial function 
$$
f(x)= \frac{x}{\delta}  \quad{(0\leq x\leq \delta)},
\quad{ f(x)=\sum_{k=0}^5 a_k x^k 
\quad{(\delta <  x \leq 2 \delta)},  }
\quad{f(x)=\frac{1-x}{1-2\delta}}
 \quad{(2 \delta< x\leq 1)}
$$
is determined by the interpolation conditions
$$
f(0)=f(1)=0, \quad{f(\delta)=f(2\delta)=1},
\quad{f'(\delta)=\frac{1}{\delta}, \quad{f'(2\delta)=- \frac{1}{1-2\delta}},}
\quad{f''(\delta)=f''(2\delta)=0}.
$$
These conditions guarantee that $ f \in   C^2[0,1] $. 
MAPLE is used to compute the coefficients:
$$  a_0 =\frac {48\,\delta -17}{2\,\delta -1} \, , \hspace*{13mm} a_1 =- \frac {183\,\delta -64}{\delta  \left( 2\,\delta -1 \right) } \, , \quad  a_2 ={\frac {12(23\,\delta -8)}{{\delta }^{2} \left( 2\,\delta -1 \right) }} \, ,  $$
$$  a_3 =-{\frac {2(99\,\delta - 34)}{{\delta }^{3} \left( 2\,\delta -1 \right) }}\, , \quad  a_4 = \frac {68\,\delta -23}{{\delta }^{4} \left( 2\,\delta -1 \right) } \, , \quad  \,\, a_5 =  -\frac {3(3\,\delta -1)}{{\delta }^{5} \left( 2\,\delta -1 \right)} \, .  $$
Straightforward computation gives
$$
\kappa (w,f) = \frac { ( 3003 +14474\,{\delta }^{3}-53525\,{\delta }^{4}-12344\,{\delta }^{ 5} ) ^{2}}{26 \, \delta   ( 858+72450\,{\delta }^{5}-531793\,{\delta }^{6 }+674178\,{\delta }^{7}) ( 2042-11999\,\delta +20182\,{\delta }^{2} )}.
$$
Since $ \delta  $ appears as a factor in the denominator, we obtain
$ \lim_{\delta \rightarrow 0} \kappa (w,f)=\infty$.  

(iii) The following example shows that the ``monotonicity'' 
requirement on $ w $ in the Corollary  
cannot be eliminated.

Let $ [a,b]=[0,1] $, $ w(x)=1-x $ and $ f(x)=\sin(\pi x/2) $. Straightforward computation gives
$$  \kappa (w,f) = \Bigl( \frac{\pi^2+4}{\pi^2 -4}\Bigr)^2 \approx  5.5835.  $$

\end{document}